\nonstopmode
\documentclass[leqno]{amsart}
\usepackage{amsthm,graphics}

\newtheorem{theorem}{Theorem}[section]

\theoremstyle{definition}
\newtheorem{defi}{Definition}[section]

\newcommand{\bc}{{\bf c}}
\newcommand{\bv}{{\bf v}}
\newcommand{\ba}{{\bf a}}
\newcommand{\bs}{{\bf s}}
\newcommand{\bw}{{\bf w}}

\newcommand{\sC}{{\mathcal C}}
\newcommand{\sP}{{\mathcal P}}
\newcommand{\T}{{\mathcal T}}
\newcommand{\NN}{{\mathbb N}}
\newcommand{\ZZ}{{\mathbb Z}}

\newcommand{\compnum}{14}      
\newcommand{\compden}{29}      
\newcommand{\comptw}{28}       
\newcommand{\comptd}{58}       
\newcommand{\compbnd}{6.14316} 
\newcommand{\comprat}{0.4827}  
\newcommand{\compdep}{60}      
\newcommand{\compones}{29}     
\newcommand{\compsnum}{28}     
\newcommand{\compsden}{59}     
\newcommand{\compsdep}{59}     
\newcommand{\compsrat}{0.4745} 

\begin{document}
\title{Lower Bounds for the Total Stopping Time of $3X + 1$ Iterates}
\author{David Applegate and Jeffrey C. Lagarias}
\address{AT\&T Laboratories, Florham Park, NJ  07932-0971}
\email[D. Applegate]{david@research.att.com}
\email[J. C. Lagarias]{jcl@research.att.com}

\date{March 7, 2001}

\begin{abstract}
The total stopping time $\sigma_{\infty}(n)$ of a positive integer
$n$ is the minimal number of iterates of the $3x+1$ function needed to
reach the value $1$, and is $+\infty$ if no iterate of $n$ reaches $1$.
It is shown that  are infinitely many positive integers 
$n$ having a finite total stopping time $\sigma_{\infty}(n)$, 
such that  $\sigma_{\infty}(n) > \compbnd \log n.$
The proof involves a search  of $3x +1$ trees to depth \compdep,
A heuristic argument suggests that for any constant 
$\gamma < \gamma_{BP} \approx 41.677647$,
a search of all $3x +1$ trees to sufficient depth could produce
a proof that there are infinitely many $n$ such that 
$\sigma_{\infty}(n) > \gamma  \log n.$   
It would require a very large computation
to search
$3x + 1$ trees to a 
sufficient depth to produce a proof that
 the expected  behavior of a ``random''
$3x +1$ iterate, which is
 $\gamma =  \frac {2}{\log 4/3} \approx 6.95212,$
occurs infinitely often. \\


\end{abstract}

\maketitle

\section{Introduction}

The  $3x+1$
problem concerns the behavior under iteration of the $3x+1$ function
$T: {\mathbf Z} \mapsto {\mathbf Z}$
given by
\begin{equation}\label{eq11}
T(n) = \left\{\begin{array}{ll}
\frac{n}{2} & \mbox{if } n \equiv 0 \pmod{2}\\
\frac{3n+1}{2} & \mbox{if } n \equiv 1 \pmod{2}
\end{array}\right.
\end{equation}
The notorious $3x+1$ Conjecture asserts that for each $n \geq 1$, some iterate
$T^{(k)} (n) = 1$. 
It has been verified for all $n < 2.702  \times 10^{16}$ by 
Oliveira e Silva~\cite{os99}, whose  subsequent computations
have extended this bound to $1.000 \times 10^{17}.$
The dynamical behavior of the $3x+1$ function 
has been extensively studied. It is a deterministic process
which nevertheless appears to exhibit a kind of pseudorandom
behavior.  Surveys on the $3x + 1$ problem can be found in 
Lagarias~\cite{l85}, M\"{u}ller~\cite{m91} and Wirsching~\cite{wi98}.

To study the behavior of iterates of
the  $3x+1$ function on the positive integers quantitatively, we make
the following definitions.
For each $n$ we call the minimal $k$ such that $T^{(k)} (n) = 1$ the
{\em total stopping time}
of $n$ and denote it $\sigma_{\infty} (n)$, letting 
$\sigma_{\infty}(n) = \infty$ if it 
is otherwise undefined. A rescaled
version of the total stopping time, which we call the {\em stopping
time ratio} $\gamma(n)$, defined by 
\begin{equation}~\label{2a} 
\gamma (n) := \frac{\sigma_{\infty}(n)}{\log n}.
\end{equation}
A trajectory {\em converges} if it has a 
finite total stopping time, and we consider it to consist
of $n$ and  all its iterates up to the first iterate 
$k$ with $T^{(k)}(n) = 1.$ Its {\em parity sequence} $\bv(n)$ is
the zero-one vector of length $\sigma_{\infty}(n)$ giving
the residue $(\bmod ~2)$ of $n$ and its iterates up to
$T^{(k-1)}(n)$.
The {\em ones-ratio} $\rho(\bv) $ of a zero-one vector $\bv$ is
the ratio of the number of ones in the vector to its length.
By extension the {\em ones ratio} $\rho(n)$ of a convergent trajectory
is defined to be the fraction of odd integers appearing
in its parity sequence, i.e. the ones-ratio $\rho(\bv(n))$. 
For example, the trajectory of $n=3$ is $(3, 5, 8, 4, 2, 1)$,
with total stopping time $\sigma_{\infty}(3) = 5$,
parity sequence $\bv(3)= (1, 1, 0, 0, 0)$ and
ones-ratio $\rho(3) = \frac{2}{5} = 0.4.$ In this case
$\gamma(3) = \frac {5}{\log 3} \approx 4.5512.$

This paper is concerned with rigorous results concerning
the behavior of the stopping time ratio  $\gamma(n).$ 
It is well known that $\gamma(n) \geq \frac{1}{\log 2}$ for all $n$, and
that equality holds exactly for $n = 2^k$ for $k \geq 1.$
In this paper our object is to get lower bounds for the
stopping time ratio,
of the form $\gamma(n) \geq \gamma,$ for some constant 
$\gamma > \frac{1}{\log 2},$
 that hold for infinitely many $n$.
We note that study of the size of the stopping time ratio  
is essentially the same as studying the size of the ones-ratio,
as given by the following bounds.
  For any convergent trajectory
one has
\begin{equation}~\label{11a} 
\gamma(n) \geq \frac{1}{\log 2 - \rho(n) \log 3}, 
\end{equation}
while if $\rho(n) \le  .61$ then,  for any positive $\epsilon$,
 $$\gamma(n) \le \frac{1}{\log 2 - \rho(n) \log 3} + \epsilon,$$
holds provided $n$ is sufficiently large, with  $n > n_0(\epsilon)$.
(It is believed that the condition
$\rho(n) \le .61$ should hold for all sufficiently large
integers, see below, while to get a uniform bound in 
terms of $\epsilon$ as above, one 
needs only assume 
$\rho(n) \le c_0 < \frac{\log 2}{\log 3} \approx .63092.$)
Our actual analysis will be based on study of the ones-ratio.

As background, we review 
the current conjectures about the
size  of the stopping time ratio.
These conjectures were suggested by 
theoretical results proved for stochastic models;
these  models make
predictions about the total stopping time function which 
successfully match experimental evidence for iterates
of the $3x + 1$ function. 
Simple models for average-case behavior of the stopping time function
appear in Crandall~\cite{c78}, Lagarias~\cite{l85}, Rawsthorne~\cite{r85},
Lagarias and Weiss~\cite{lw92}
and Borovkov and Pfeifer~\cite{bp93}. These models predict that
the average value of $\gamma(n)$ should
be 
$$\frac {2}{\log 4/3} \approx 6.95212,$$
which corresponds to a ones-ratio of $\frac{1}{2},$ and that
nearly all  values
of $\gamma(n)$ should be within 
$O( (\frac {1}{\sqrt{\log n}} )^{1 + \epsilon})$
of this value. This is well supported by experimental evidence.

Lagarias and Weiss \cite{lw92} also developed  stochastic models to describe
the maximal values attained by the total stopping time. These models are
are analyzed using the theory of  large deviations. They  predict
that  the limit
superior of $\gamma (n)$ as $n \to \infty$ should
be a certain constant $\gamma_{BP} \approx 41.677647$,
which corresponds to a ones-ratio of about $0.609091.$
(This constant is the solution to a certain
functional equation given in \cite{lw92}.) The occurrence of values close
to the extremal one are extremely rare events whose probability can
be estimated in the stochastic models, see
\cite[Section 6]{lw92}). For example, the largest value
to be expected for $n < 10^{17}$ is about $32,$  and
Oliveira e Silva's record number is
$n = 1, 008, 932, 249, 296, 231$
with $\sigma_\infty(n) = 1142 $ and $\gamma(n) \approx 33.0558.$
Backwards search methods have uncovered
some larger integers which have larger values of $\gamma(n)$.
V. Vyssotsky~\cite{v97} found that $\gamma(n) \approx 35.2789$ for 
$$n = 37, 664, 971, 860, 959, 140, 595, 765, 286, 740, 059,$$
which has $\sigma_{\infty}(n) = 2565$. He also found the 
 current record value for $\gamma(n)$ of about $36.40$, 
coming from a number around $10^{110}$.

We now turn to the problem considered here, that of obtaining
bounds $\alpha$ such that there are provably infinitely many
positive $n$ having $\gamma(n) \geq \gamma.$ 
The current best
value of $\alpha$ is that associated to the family of numbers
$n = 2^k - 1.$ It is well known that after $k$ iterations one
reaches $T^{(k)}(n) = 3^k - 1.$ Since this number cannot get to
$1$ any faster than a power of $2$ (which in fact occurs for $k=2$),
one obtains the rigorous lower bound
\begin{equation}~\label{eq12}
 \gamma (2^k - 1) \geq \frac { \log 2 + \log 3}{(\log 2)^2} \approx 3.729.
\end{equation}
This corresponds to a ones ratio about $0.387$. 
It is expected that the numbers $3^k - 1$ for most $k$ have  iterates
behaving like a ``random'' integer, in which case the lower
bound (\ref{eq12}) could be significantly improved for infinitely many
$k$.
However, as  far as is currently known,  it could be the case that for
each $\epsilon > 0$ the bound
$$\sigma_{\infty} ( 3^k - 1) < (\frac{1}{\log 2} + \epsilon) \log (3^k - 1),$$
holds for all sufficiently large $k$; if so 
then (\ref{eq12}) could not
be improved upon, asymptotically as $k \to \infty.$. 

The object of this paper is to improve this lower bound, as follows.

\begin{theorem}\label{th11}
There  exist infinitely many $n \ge 1$ having a 
convergent trajectory with a ones-ratio of 
at least $\frac{\compnum}{\compden}\approx  \comprat$.
Thus there  is an infinite set of positive integers $n$ such that
$\gamma(n)$ is finite and
\begin{equation}\label{eq13}
\gamma(n) \geq  \frac{\compden}{\compden\log2 - \compnum \log 3} \approx \compbnd.
\end{equation}
\end{theorem}

This result is proved by a study of  all $3x+1$ trees to depth $\compdep$,
by extensive computation, involving paths with up to $\compones$ ones.
In the rest of the paper we describe $3x+1$ trees, the form of proof
certificates for lower bounds of the kind in Theorem~\ref{th11} and
the computations.

%
%

\section{3X + 1 Trees}

The $3x+1$ iteration when run backwards from any fixed integer $a$
produces a tree of preimages of $a$.  

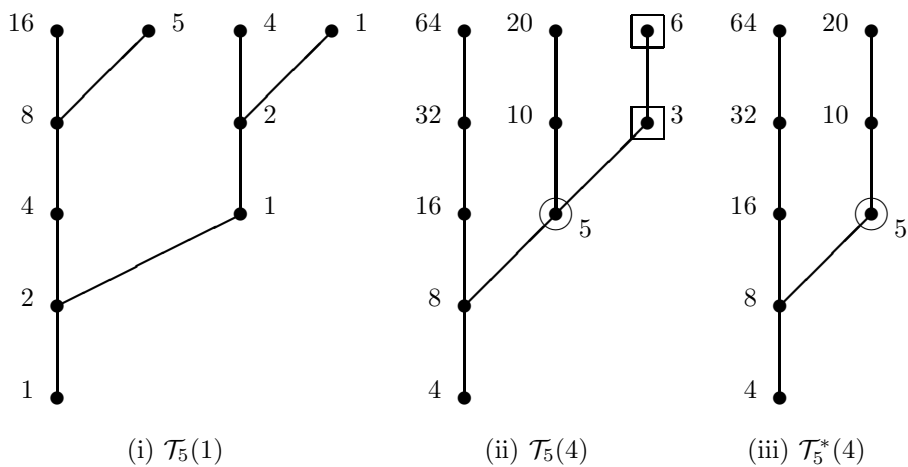
\begin{figure}[hbp]\centering
\setlength{\tabcolsep}{0.125in}
\vspace{.1in}
\begin{tabular}{ccc}
\setlength{\unitlength}{0.012in}
\begin{picture}(153,172)(57,436)
\thicklines
\put( 70,600){\makebox(0,0)[rb]{16}}
\put(210,600){\makebox(0,0)[lb]{1}}
\put(170,520){\makebox(0,0)[lb]{1}}
\put(170,560){\makebox(0,0)[lb]{2}}
\put(170,600){\makebox(0,0)[lb]{4}}
\put(160,560){\line( 1, 1){ 40}}
\put(160,560){\line( 0, 1){ 40}}
\put( 80,480){\line( 2, 1){ 80}}
\put(160,520){\line( 0, 1){ 40}}
\put( 70,520){\makebox(0,0)[rb]{4}}
\put(130,600){\makebox(0,0)[lb]{5}}
\put( 70,560){\makebox(0,0)[rb]{8}}
\put( 70,480){\makebox(0,0)[rb]{2}}
\put( 70,440){\makebox(0,0)[rb]{1}}
\put( 80,440){\circle*{6}}
\put(120,600){\circle*{6}}
\put( 80,600){\circle*{6}}
\put( 80,560){\circle*{6}}
\put( 80,520){\circle*{6}}
\put( 80,480){\circle*{6}}
\put( 80,560){\line( 1, 1){ 40}}
\put(160,520){\circle*{6}}
\put( 80,560){\line( 0, 1){ 40}}
\put( 80,520){\line( 0, 1){ 40}}
\put( 80,480){\line( 0, 1){ 40}}
\put( 80,440){\line( 0, 1){ 40}}
\put(200,600){\circle*{6}}
\put(160,600){\circle*{6}}
\put(160,560){\circle*{6}}
\end{picture}
&
\setlength{\unitlength}{0.012in}
\begin{picture}(113,172)(57,436)
\thicklines
\put(120,520){\line( 1, 1){ 40}}
\thinlines
\put(167,593){\line( 0, 1){ 14}}
\put(167,607){\line(-1, 0){ 14}}
\put(153,607){\line( 0,-1){ 14}}
\put(153,593){\line( 1, 0){ 14}}
\put(167,553){\line( 0, 1){ 14}}
\put(167,567){\line(-1, 0){ 14}}
\put(153,567){\line( 0,-1){ 14}}
\put(153,553){\line( 1, 0){ 14}}
\thicklines
\put(120,560){\line( 0, 1){ 40}}
\put( 80,480){\line( 1, 1){ 40}}
\put( 70,600){\makebox(0,0)[rb]{64}}
\put(130,510){\makebox(0,0)[lb]{5}}
\put(110,560){\makebox(0,0)[rb]{10}}
\put(110,600){\makebox(0,0)[rb]{20}}
\put(170,600){\makebox(0,0)[lb]{6}}
\put(170,560){\makebox(0,0)[lb]{3}}
\put(120,520){\line( 0, 1){ 40}}
\thinlines
\put(120,520){\circle{14}}
\thicklines
\put(120,600){\circle*{6}}
\put( 70,560){\makebox(0,0)[rb]{32}}
\put( 70,520){\makebox(0,0)[rb]{16}}
\put( 70,480){\makebox(0,0)[rb]{8}}
\put( 80,440){\circle*{6}}
\put( 80,480){\circle*{6}}
\put( 80,520){\circle*{6}}
\put( 80,560){\circle*{6}}
\put( 80,600){\circle*{6}}
\put( 70,440){\makebox(0,0)[rb]{4}}
\put( 80,440){\line( 0, 1){ 40}}
\put( 80,560){\line( 0, 1){ 40}}
\put( 80,520){\line( 0, 1){ 40}}
\put( 80,480){\line( 0, 1){ 40}}
\put(160,560){\circle*{6}}
\put(160,560){\line( 0, 1){ 40}}
\put(160,600){\circle*{6}}
\put(120,560){\circle*{6}}
\put(120,520){\circle*{6}}
\end{picture}
&
\setlength{\unitlength}{0.012in}
\begin{picture}(73,172)(57,436)
\thinlines
\put(120,520){\circle{14}}
\thicklines
\put(120,560){\line( 0, 1){ 40}}
\put(120,520){\line( 0, 1){ 40}}
\put( 80,480){\line( 1, 1){ 40}}
\put(130,510){\makebox(0,0)[lb]{5}}
\put(110,600){\makebox(0,0)[rb]{20}}
\put(110,560){\makebox(0,0)[rb]{10}}
\put( 70,600){\makebox(0,0)[rb]{64}}
\put( 70,480){\makebox(0,0)[rb]{8}}
\put( 70,520){\makebox(0,0)[rb]{16}}
\put( 70,440){\makebox(0,0)[rb]{4}}
\put( 70,560){\makebox(0,0)[rb]{32}}
\put( 80,440){\circle*{6}}
\put( 80,600){\circle*{6}}
\put( 80,560){\circle*{6}}
\put( 80,520){\circle*{6}}
\put( 80,480){\circle*{6}}
\put( 80,560){\line( 0, 1){ 40}}
\put(120,600){\circle*{6}}
\put( 80,520){\line( 0, 1){ 40}}
\put( 80,480){\line( 0, 1){ 40}}
\put( 80,440){\line( 0, 1){ 40}}
\put(120,520){\circle*{6}}
\put(120,560){\circle*{6}}
\end{picture}
\\
(i) $\T_5(1)$ & (ii) $\T_5(4)$ & (iii) $\T^*_5(4)$
\end{tabular}
\caption{$3x+1$ trees $\T_k(a)$ and ``pruned'' $3x+1$ tree $\T^*_k(a)$}
\label{fig1}
\end{figure}

A {\em $3x+1$ tree} $\T_k(a)$ is a rooted, labelled
tree of depth $k$, representing the inverse iterates $T^{-j}(a)$ for
$0 \leq j \leq k$.  The inverse map $T^{-1}(n)$ is multivalued:
$$
T^{(-1)}(n) = \left\{\begin{array}{ll}
\{2n\} & \mbox{if $n \equiv 0$ or $1 \pmod{3}$,}
\vspace{0.1in}\\
\{2n,\frac{2n-1}{3}\} & \mbox{if $n \equiv 2 \pmod{3}$.}
\end{array}\right.
$$
The root node $a$ is at depth 0, and a node labelled $n$ at depth $d$
of the tree is connected by an edge to a node labelled $T(n)$ at depth
$d-1$ of the tree\footnote{We adopt a convention of ``unrolling'' any
cycles under $T$, so that the same node label may appear at different
levels of the tree if a cycle is present, cf. Figure~\ref{fig1}.}. 
Thus the {\em depth}  of a node is the number of edges in  a path from  
it to the root node.
As described in~\cite{l85}, the nodes labelled $n
\equiv 0 \pmod{3}$ give rise only to a linear chain of nodes labelled
$n' \equiv 0 \pmod{3}$ at higher levels.  It is convenient to remove
all such nodes and study a ``pruned'' tree $\T^{*}_{k}(n)$ consisting
of nodes $n \not\equiv 0 \pmod{3}$.  That is, we study the
inverse map on the set of positive integers $n \not\equiv 0 \pmod{3}$
given by
$$
{T^*}^{(-1)}(n) = \left\{\begin{array}{ll}
\{2n\} & \mbox{if $n \equiv 1, 4, 5$ or~ $7  \pmod{9}$,}
\vspace{0.1in}\\
\{2n,\frac{2n-1}{3}\} & \mbox{if $n \equiv 2$ or $8 \pmod{9}$.}
\end{array}\right.
$$
Figure~\ref{fig1} presents some
examples of $\T_k(a)$ and $\T^*_k(a)$.  (Nodes $n \equiv 5 \pmod{9}$
are circled to indicate that they have some preimage $T^{-1}(n) \equiv
0 \pmod{3}$, and nodes $n \equiv 0 \pmod{3}$ are indicated with a square.)

We assign to each edge of a $3x+1$ tree between depth $d$ and
$d - 1$ an {\em edge label}
$0$ or $1$, which is the value $n~ (\bmod~2)$ of the endpoint
node at depth $d$. Traversing a
path from a leaf node labelled $n$ in $\T_k^*(a)$ at depth $d$
to the root node $a$, the successive edge labels encode the 
initial part of the parity
vector $\bv(n)$ of $n$, and the zero-one weight vector $\bw(n; a)$
is the reversal of the vector $\bv(n)$ to the root and gives the
sequence of edge labels from the root node to the node $n$.
The {\em weight} $w$ of such a path is the
sum of the edge labels on the path, and we sometimes say that
the node $n$ has {\em weight} $w(n)$.
Note that $\frac{w(n)}{d}$ is the
ones-ratio of this vector along the  path. 
The {\em max-weight} $w_{max}(a)$
of  a tree $\T_k^*(a)$ is the largest weight among all paths 
from a leaf node to the root node. 
Thus $\frac{w_{max}(a)}{k}$ is the maximal 
ones-ratio among any path in the tree  $\T_k^*(a)$ from a
depth $k$ node to the root.

The {\em structure}  of a pruned $3x+1$ tree $\T_k^*(a)$
of depth $k$ is its isomorphism class as a
rooted, edge-labelled tree. That is, we
 say that two pruned $3x+1$ trees $\T_k^*(a)$ and $\T_k^*(b)$ have the
same {\em structure} if they are isomorphic as rooted trees by an
isomorphism that preserves edge labels.
The structure of a 
``pruned'' $3x+1$ tree
$\T_k^*(a)$ is completely determined by knowledge of $a~(\bmod~ 3^{l + 2})$,
where $l = w_{max}(a).$ 
To demonstrate this, we grow the tree from the
root, labelling each new node $n~(\bmod~3^m)$, where $m$ is determined
recursively, as follows. The root node is assigned $m = w_{max}(a) + 2.$
For each
node $n ~(\bmod~3^m)$ at level $d$, we assert $m \ge 2$,
and if $n \equiv 1, 4, 5$ or $7 ~(\bmod~9)$ there is a single
edge labelled $0$ to a node at level $d+1$ labelled $2n ~(\bmod~3^m)$.
 If $n \equiv 2 $ or $8 ~(\bmod~9)$, there is one node
at depth $d+1$ as above plus a second edge labelled $1$
to a node at depth $d+1$ labelled $ (2n - 1)/3 ~(\bmod~3^{m - 1})$.
Since we started with $m = w_{max}(a) + 2$ at the root node, and
since no path has more than $w_{max}(a)$ edge labels of $1$, the
value of the exponents $m$ stays at least $2$ all the way to,
and including, the leaves at level $k$.

We can improve on this bound, in
one circumstance. We  say that a tree $\T_k^*(a)$ is
{\em critical} at max-weight (or max-level) $l$ if it
has max-weight $w_{max}(a)= l$  and if all nodes at depth $k - 1$
have weights $w_{max}(a) - 1$ or less in their path to the root.
The structure of a critical tree is completely determined
by $a~( \bmod~ 3^{l + 1})$ where $l= w_{max}(a)$,
rather than  $a~(\bmod~ 3^{l + 2})$. The proof is the 
same as before, noting that the exponent $m \geq 2$ on every node up to
and including depth $k - 1$, but $m = 1$ may occur on some leaf
nodes, on paths of weight exactly $w_{max}(a)$. 

For trees of depth $k$ one has $l \le k$, and the
case $l = k$ is necessarily critical, so  by 
the  discussion above all pruned tree structures  $\T_k^*(a)$
are determined by $a~(\bmod~ 3^{k + 1})$. Thus 
 there are 
at most $2 \cdot 3^{k}$ distinct pruned tree structures $\T_k^*(a)$
of depth $k$. The actual number $R(k)$ of distinct tree structures 
is smaller but still grows exponentially in $k$.

We now suppose given as data a
congruence class $a ~(\bmod~ 3^{l + 1})$, and from it we may grow a 
``pruned'' $3x+1$ tree with root $a$ to the (unique) depth $k$
at which it becomes critical with max-weight $l$. We represent
$a ~(\bmod~ 3^{l + 1})$ as a ternary (or $3$-adic) vector
$\bc = \bc(a) = (c_{0}, c_{1}..., c_l)$ of length $l + 1$, where
$$ a = (c_{l}c_{l-1}... c_0)_3 = \sum_{j=0}^{l} c_j 3^j, $$
with $c_j \in \{0, 1, 2\}.$ 
We label
this tree $\T_{[l]}^*(\bc)$. Using the procedure above, we can
construct node labels for this tree consisting
of such $3$-adic vectors of variable length $m$ at each node. 
Figure~\ref{fig2} illustrates this,  depicting
the ``pruned'' tree $\T_{5}^*(4)$, which is critical with  max-weight $2$,
and the same tree $\T_{[2]}(011)$ computed using root value 
$4~(\bmod~27).$ 
\begin{figure}[htb]
\begin{center}
\input depth.pstex_t
\end{center}

\caption{``Pruned'' $3x+1$ tree $\T_{[k]}^\ast (\bc )$ with $\bc= 110$,
for $a = (011)_3$.}
\label{fig2}
\end{figure}

We define the {\em ones-ratio} $\rho( \bc)$ for a 
critical tree $\T_{[l]}^*(\bc)$ to
be 
$$\rho( \bc) : = \frac {l}{k}, $$
where $k$ is the depth of the tree, and $l$ is its criticality level.

%
%

\section{Lower Bound Certificates}

We describe a finite set of information (``certificates'') 
that are sufficient to yield proof that there are infinitely
many $n \geq 1$ having a finite total stopping time whose
ones-ratio satisfies the lower bound
\begin{equation}\label{Eq31}
\rho (n) \ge \alpha ~.
\end{equation}

\begin{defi}\label{de31}
A {\em prefix code} $\sC$ for ternary sequences is a 
finite set of
codewords $\bc = (c_0, c_1, \ldots, c_l ) \in \{0,1,2\}^{l+1}$
 of varying lengths
$l= l(\bc)$ having the properties:
\begin{itemize}
\item[(i)]
If $\bc, \bc' \in \sC$ with $\bc \neq \bc'$, then $\bc$ 
is not a prefix of $\bc'$ or vice-versa.
\item[(ii)]
Any infinite sequence 
$\bs = (s_1, s_2, s_3 , \ldots ) \in \{0,1,2 \}^{\NN_0}$ 
has an initial segment $(s_1, s_2, \ldots, s_{l+1} )$
that is a codeword in $\sC$.
\end{itemize}
\end{defi}

\begin{defi}\label{de32}
A ({\em lower bound) certificate for ones-ratio} $\alpha$ consists of:
\begin{itemize}
\item[(i)]
A finite list $\sC$ of residue classes $a~(\bmod~3^{l+1} )$ for 
various lengths,
$l= l(a)$, which are exhaustive, in the
sense that  their associated vectors
$\bc (a)$, together with the additional vector $\bc_0 = 0$,
 form a (ternary)  prefix code.
\item[(ii)]
With each $a \in \sC$ with $a~(\bmod~3^{l+1} )$ is given a 
{\em critical path vector}
$\bw(a) \in \{0,1\}^k$, where $k$ satisfies
$$
l/k \ge \alpha ~,
$$
such that $\bw (a)$ has weight at most $l$, and if it has weight 
exactly $l$, then a one occurs in its last coordinate, and 
$\bw (a)= \bw(n; a)$ for some leaf node of this $3x+1$ tree, i.e.
$\bw(a)$ is  a labelled path in the $(3x+1 )$ tree of 
$a~(\bmod~3^{l+1} )$, from the root node to some  node at depth $k$.
\end{itemize}
The certificate vector certifies that the associated critical tree
$\T_{[l]}^\ast (\bc( a ))$ with max-level $l$ has depth\footnote{A
single path is a $3x+1$ tree sometimes can be grown to a depth 
greater than $k$,
from knowledge of $a~(\bmod ~3^{l+1} )$, as long as 
it has weight at most $l$, and if exactly weight $l$ then it 
has a one in its last coordinate.}
$k'$ with $k' \le k$ and has ones-ratio $l/k' \ge \alpha$.
\end{defi}

As an example, we present in  Table  \ref{tab1}
a  certificate for level $\alpha = 1/3$. 
\begin{table}
\center
\renewcommand{\arraystretch}{0.625}
\begin{tabular}{|r|r|r|l|}
\hline
\rule[-0.07in]{0cm}{.23in}
$a \;(\mbox{mod}~ 3^{l+1})$ & $l$ & $k$ & critical path\\
\hline
   01 & 1 & 2 & 01 \\
   11 & 1 & 2 & 01 \\
  021 & 2 & 6 & 000101 \\
  121 & 2 & 6 & 000101 \\
 0221 & 3 & 9 & 000001101 \\
 1221 & 3 & 9 & 000001101 \\
02221 & 4 &12 & 000100000111 \\
12221 & 4 &11 & 00010001011 \\
22221 & 4 &11 & 00010001011 \\
   02 & 1 & 1 & 1 \\
   12 & 1 & 3 & 001 \\
   22 & 1 & 1 & 1 \\
\hline
\end{tabular}
\caption{Certificate for $\rho = \frac{1}{3}$.}
\label{tab1}
\end{table}
In this table the numbers $a~(\bmod~3^{l+1})$ are represented
in their ternary expansion, so form a suffix code, rather than
a prefix code.
(The vectors $\bc (a)$ reverse the order of the digits.)
This certificate has max-weight $4$ and max-depth $12$,
and consists of $12$ different trees.

\begin{defi}\label{de33}
A {\em strong lower bound certificate for ones-ratio} $\alpha$ consists of a 
finite list $\sC$ as above satisfying:
\begin{itemize}
\item[(i')] Same as (i) above.
\item[(ii')]
With each $a~(\bmod~3^{l+1} ) \in \sC$ are given two 
critical path  vectors $\bw_j (a) \in \{0,1\}^{k_j}$, with
weights $l_j \le l$, respectively, such that neither is a prefix of
the other, and where $k_j$ and $l_j$ satisfy
$$
l_j/k_j \ge \alpha ~,
$$
and if $\bw_j$ has weight exactly $l$
then it has a 1 in its last coordinate, and each $\bw_j (a)$ 
can be found as a labeled path of length $k_j$ from the root
in the $3x+1$ tree of $a~(\bmod~3^{l + 1} )$.
\end{itemize}
\end{defi}

As an example, there  exists a strong lower bound
certificate for $\alpha = \frac{1}{3}$
which has max-weight $4$ and max-depth $12$, the same as the
lower bound certificate above, but which 
consists of $42$ trees. It appears in Table~\ref{tab2} below.

\begin{table}
\renewcommand{\arraystretch}{0.625}
\begin{tabular}{|rr|rl|rl|}
\hline
\rule[-0.07in]{0cm}{.23in}
$a \;(\mbox{mod}~ 3^{l+1})$ & $l$ & $k_1$ & critical path 1 & $k_2$ &
critical path 2\\
\hline
   001 & 2 & 4 & 0101 & 6 & 010001 \\
   101 & 2 & 4 & 0101 & 5 & 00011 \\
   201 & 2 & 5 & 00011 & 6 & 010001 \\
  0011 & 3 & 7 & 0100101 & 8 & 01000011 \\
  1011 & 3 & 7 & 0100101 & 9 & 010010001 \\
  2011 & 3 & 8 & 01000011 & 9 & 010010001 \\
   111 & 2 & 3 & 010 & 3 & 011 \\
   211 & 2 & 3 & 011 & 5 & 01001 \\
  0021 & 3 & 8 & 00010101 & 9 & 000100011 \\
  1021 & 3 & 6 & 000101 & 9 & 000100011 \\
 02021 & 4 & 10 & 0001010101 & 12 & 000101010001 \\
 12021 & 4 & 10 & 0001010101 & 11 & 00010100011 \\
 22021 & 4 & 11 & 00010100011 & 12 & 000101010001 \\
  0121 & 3 & 8 & 00000111 & 9 & 000101001 \\
  1121 & 3 & 7 & 0001011 & 8 & 00000111 \\
  2121 & 3 & 7 & 0001011 & 9 & 000101001 \\
 00221 & 4 & 11 & 00000100111 & 12 & 000001101001 \\
 10221 & 4 & 10 & 0000011011 & 12 & 000001101001 \\
 20221 & 4 & 10 & 0000011011 & 11 & 00000100111 \\
 01221 & 4 & 11 & 00000110101 & 12 & 000001100011 \\
 11221 & 4 & 11 & 00000110101 & 12 & 000100010101 \\
 21221 & 4 & 9 & 000001101 & 12 & 000100010101 \\
002221 & 5 & 14 & 00010000011101 & 15 & 000100000110011 \\
102221 & 5 & 12 & 000100000111 & 15 & 000100010100101 \\
202221 & 5 & 14 & 00010000011101 & 15 & 000100010100101 \\
 12221 & 4 & 11 & 00010001011 & 12 & 000100000111 \\
 22221 & 4 & 11 & 00010001011 & 12 & 000001001101 \\
    02 & 1 & 1 & 1 & 3 & 001 \\
   012 & 2 & 5 & 00101 & 6 & 000011 \\
  0112 & 3 & 7 & 0010101 & 8 & 00100011 \\
  1112 & 3 & 7 & 0010101 & 9 & 001010001 \\
  2112 & 3 & 8 & 00100011 & 9 & 001010001 \\
   212 & 2 & 3 & 001 & 6 & 000011 \\
   022 & 2 & 4 & 1001 & 6 & 100001 \\
   122 & 2 & 2 & 11 & 4 & 1001 \\
   222 & 2 & 2 & 11 & 6 & 100001 \\
\hline
\end{tabular}
\caption{Strong certificate for $\rho = \frac{1}{3}$}
\label{tab2}
\end{table}

\begin{theorem}\label{th31}
(i) The existence of a lower bound certificate for ratio $\alpha$ implies:
Given any $a \ge 1$ with $a \not\equiv 0$ $(\bmod~3)$, 
not in a cycle, there exist infinitely many $n$ having
\begin{equation}\label{Eq35}
T^{(k)} (n) =a ~,
\end{equation}
for some $k \ge 1$, and such that
 the ones-ratio $\rho (n; a) := \rho(\bw(n; a))$ 
of the path from $n$ to $a$ has
\begin{equation}\label{Eq36}
\rho (n; a) = l/k \ge \alpha ~,
\end{equation}
(ii) The existence of a strong lower bound certificate for 
ones-ratio $\alpha$ implies: Given any $a \ge 1$ with 
$a \not\equiv 0$ $(\bmod~3)$ the number of
$n \le x$ with $T^{(k)} (n) =a$ for some $k \ge 1$ and with
\begin{equation}\label{Eq37}
\rho (n; a) = l/k \ge \alpha ~,
\end{equation}
is at least $c(a) x^\delta$ for certain $c(a)>0$
and  $\delta = \delta(a) > 0 $. 
\end{theorem}

\begin{proof}
(i) Since $a$ is not in a cycle, all the nodes in the (infinite) 
``pruned'' $3x +1$ tree $\T^\ast (n)$ with root node $a$ have distinct values.

We can find an infinite path
 $\bs = (s_1, s_2, s_3, \ldots ) \in \{0,1\}^{\NN_0}$ from 
the root node $a$ in $\T^\ast (a)$ which contains infinitely many nodes
satisfying (\ref{Eq36}).
View $a$ as a 3-adic integer, encoded as
$$\bc (a) = (c_0,c_1, c_2, c_3, \ldots ) \in \{ 0,1,2 \}^{\NN_0} ~;$$
in which $c_0 \neq 0$ and
$c_j =0$ for all sufficiently large $j$.
Then $\bc (a)$ has a prefix in the certificate $\sC$, of length $l+1$, say.
By property (ii) of the certificate, there exists a path in the associated
tree $\T_{[l_1]}^\ast (\bc (a))$ from leaf node $n$ to $a$ with
$$T^{(k_1)} (n_1) = a$$
and
\begin{equation}\label{Eq38}
\rho (n; a ) = \frac{l_1}{k_1} \ge \alpha ~,
\end{equation}
where $k_1$ is the depth of $\T_{[l]}^\ast (c ( \ba ))$.
Now $n_1 \in \ZZ_{\ge 0}$ is not in a cycle and 
$n_1 \not\equiv 0$ $(\bmod~3)$, so we can repeat the same construction;
finding a prefix of $\bc (n_1 )$ of length $l_2 +1$, say, with tree
$\T_{[l_2]}^\ast (\bc (n_1 ))$ of depth $k_2$ having a path from leaf node
$n_2$ to $n_1$ with
$$T^{(k_2)} (n_2) = n_1$$
and
$$\rho (n_2, n_1) = \frac{l_2}{k_2} \ge \alpha \,.$$
Thus
$$T^{(k_1 + k_2 )} (n_2) = a$$
and
\begin{equation}\label{Eq39}
\rho (n_2, a ) \ge \min \left[
\frac{l_2}{k_2} , \frac{l_1}{k_1} \right] \ge \alpha ~.
\end{equation}
Continuing
inductively, we find an infinite chain 
$a \leftarrow n_1 \leftarrow n_2 \leftarrow n_3 \leftarrow \ldots$ 
for which (\ref{Eq35}) and (\ref{Eq36}) hold.

(ii) We use a similar argument.
Suppose first that $a$ is not in a cycle.
Then at the first round we find two different paths in 
$\T_{[l_1]}^\ast (\bc (a))$ with $n_1 \to a$, $n'_1 \to a$ 
satisfying (\ref{Eq36}).
At the $j$-th stage we produce $2^j$ distinct elements 
$n_j \to a$ satisfying (\ref{Eq36}).
Furthermore, if the largest depth tree $\T_{[l]}^\ast (\bc )$ 
for $\bc \in \sC$ is $d$ then, all the elements $n_j$ are at 
depth at most $jd$ from the root node $a$, hence 
each such element satisifes $n_j \leq 2^{jd} a$.
Let $\pi_a^{\alpha} (x)$ count the number of elements $\le x$ which satisfy
(\ref{Eq35}) and (\ref{Eq36}).
It follows that for $2^{jd} a \le x < 2^{(j+i)d} a$ we have
\begin{equation}\label{Eq310}
\pi_a^\alpha (x) \ge 2^j > \frac{1}{2} \left( \frac{x}{a} \right)^{1/d} = 
c(a) x^{1/d}
\end{equation}
with $c(a) = \frac{1}{2} a^{-1/d}$.

Now suppose that $a$ is in a cycle, i.e. is periodic.
Then in the first tree $\T_{[l]}^\ast ( \bc (a))$ at least one of
the two preimages $n_1$ and $n'_1$ of $a$ that satisfies (\ref{Eq37})
cannot be in the periodic cycle containing $a$.
Let it be $n_1$, and we may then apply the argument 
above to the tree with root $n_1$ for which
$$\pi_{n_1}^\alpha (x) \ge \frac{1}{2} (n_1)^{-1/d} x^{1/d} ~.$$
But any such element $n \to n_1 \to a$ by adjoining the 
path from $n_1$ to $a$,
and the inequality (\ref{Eq39}) gives
$$\rho (n,a) \ge \min
[\rho (n,n_1 ), \rho (n_1,a)] \ge \alpha ~.
$$
Thus
\begin{equation}~\label{Eq311}
\pi_a^\alpha (x) \ge \pi_{n_1}^\alpha (x) \ge \frac{1}{2} (n_1)^{-1/d}
x^{1/d} \ge \frac{1}{4} a^{-1/d} x^{1/d} ~,
\end{equation}
since $n_1 \le 2^d a$.
\end{proof}

Theorem \ref{th31} reduces the problem of finding lower bounds for 
$\rho (n)$ for infinitely many $n$ to that of finding suitable certificates, 
resp. strong certificates.
Such certificates can be searched for by computer, as we describe next.
We note that for any fixed 
$\alpha > \frac{1}{2}$ the stochastic models studied
in Lagarias and Weiss ~\cite{lw92} predict that the number of $n \le x$
that have $\rho (n) \geq \alpha$ is bounded above by 
$x^{\gamma(\alpha) + o(1)}$
for a certain exponent $\gamma(\alpha) < 1.$ 
(The value  $\gamma(\alpha)$ is effectively
computable as the fixed point of a certain functional equation.)
It follows that
strong lower bound certificates for such $\alpha$, if they could
be found, would produce a lower bound
 qualitatively of the correct form as what 
would be expected from the stochastic model predications, but 
(presumably) with
a much smaller exponent than $\gamma(\alpha)$.

%
%

\section{Computational Results}

We may search for lower bound certificates resp. strong lower bound
certificates by a ``greedy'' algorithm, as follows.

\subsection*{Certificate Search Algorithm}
\begin{itemize}
\item[(0)]
{\em Input Data}:
$\alpha$ with $0 < \alpha < 1$.
\item[(1)]
{\em Initialization}:
$\sC = \{ \bc_{ij} : i \in \{1,2\} , j \in \{0,1,2\}\}$
with $\bc_{ij} = (i,j) \in \{0,1,2\}^2$.
Set level $l=1$.
Declare all vectors $\bc_{ij} \in \sC$ {\em open}.
\item[(2)]
{\em Tree Test Step}:
Select an open vector $\bc \in \sC$, 
and let $l+1$ be the length of $\bc$.  
Compute $\T_l^\ast ( \bc )$ and determine
its depth $d= d(l, \bc )$.
If $l/d \ge \alpha$ declare the vector $\bc$ {\em closed}.
Otherwise split $\bc$ into three vectors of length $l+2$, 
$(\bc, 0)$, $(\bc, 1)$ and $(\bc, 2)$, at level $l+1$.
Update $\sC$ by deleting $\bc$ and adding the three new vectors.
Declare the three new vectors {\em open}.
\item[(3)]
{\em Termination Test}:
If all $\bc \in \sC$ are closed, halt.
Otherwise, repeat the tree test step.
\end{itemize}

In order to keep the set
of open prefixes $\sC$ as small as possible during this algorithm,
it is convenient to choose
at each test step to examine an open prefix vector that is
one of those currently of greatest level.

This algorithm is not guaranteed to halt, but if it does, it produces
a lower bound certificate for ratio $\alpha$.
The correctness of the algorithm is based on the observation 
that $\sC$ is a prefix code at all times:
the ``splitting'' step replacing $\bc$ by $(\bc,0)$, $(\bc, 1)$ 
and $(\bc,2 )$ preserves the prefix code property.
It is also easy to see that the algorithm will find a certificate 
$\sC$ if any certificate exists for ones-ratio $\alpha$,
and it will be one of minimal depth.

The algorithm above is easily modified to give a 
{\em Strong Certificate Search Algorithm}.
It is exactly the same, except that the rule for closing a node is modified:
A vector $\bc \in \sC$ at level $l$ is declared
{\em closed} (only) if $l/d \ge \alpha$ and at least 
two paths of weight $l$ exist from the root to a leaf of $\T_l^\ast (\bc )$.
If the algorithm halts, it produces a strong lower 
bound certificate for ones-ratio $\alpha$.

We searched by computer for certificates and strong certificates,  
for different values of $\alpha$.
The  certificate and strong certificate for level 
$\alpha = 1/3$ presented in Table \ref{tab1}
and Table ~\ref{tab2} were found in this way. 
The size of certificates grows rapidly with increasing 
$\alpha$.
Table~\ref{tab3} presents data giving the maximal value of $\alpha$ 
for which a certificate exists in which all trees have a critical
path containing at most $l$ ones, for $ 1 \le l < \compones$, as well
as a (not necessarily maximal) value of $\alpha$ for which a
certificate exists for $l = \compones$.

\begin{table}
\renewcommand{\arraystretch}{0.625}
\begin{tabular}{|r|r|r|r|}
\hline
\rule[-0.07in]{0cm}{.23in}
ones-ratio $\alpha$ & size $|\sC|$ & max-weight $l$ & max-depth $k$ \\
\hline
$1/4   =       0.250 $ &         6 &  1 &  4 \\
$2/7   \approx 0.286 $ &         8 &  2 &  7 \\
$3/10  =       0.300 $ &        10 &  3 & 10 \\
$1/3   \approx 0.333 $ &        12 &  4 & 12 \\
$5/14  \approx 0.357 $ &        34 &  5 & 14 \\
$5/14  \approx 0.357 $ &        34 &  6 & 14 \\
$7/19  \approx 0.368 $ &        68 &  7 & 19 \\
$8/21  \approx 0.381 $ &       120 &  8 & 21 \\
$9/23  \approx 0.391 $ &       268 &  9 & 23 \\
$2/5   =       0.400 $ &       276 & 10 & 25 \\
$11/27 \approx 0.408 $ &       704 & 11 & 27 \\
$12/29 \approx 0.414 $ &      1522 & 12 & 29 \\
$13/31 \approx 0.419 $ &      2404 & 13 & 31 \\
$14/33 \approx 0.424 $ &      4758 & 14 & 33 \\
$3/7   \approx 0.429 $ &      4782 & 15 & 35 \\
$16/37 \approx 0.432 $ &     10646 & 16 & 37 \\
$15/34 \approx 0.441 $ &     42336 & 17 & 38 \\
$4/9   \approx 0.444 $ &     48718 & 18 & 40 \\
$19/42 \approx 0.452 $ &    282326 & 19 & 42 \\
$5/11  \approx 0.455 $ &    285098 & 20 & 44 \\
$21/46 \approx 0.457 $ &    519802 & 21 & 46 \\
$11/24 \approx 0.458 $ &    829044 & 22 & 48 \\
$23/50 =       0.460 $ &   1413986 & 23 & 50 \\
$7/15  \approx 0.467 $ &   4303530 & 24 & 51 \\
$8/17  \approx 0.471 $ &   8035246 & 25 & 53 \\
$26/55 \approx 0.473 $ &  16669294 & 26 & 55 \\
$9/19  \approx 0.474 $ &  16671812 & 27 & 57 \\
$28/59 \approx 0.475 $ &  26948336 & 28 & 59 \\
$\geq 14/29 \approx 0.483 $ & 350688758 & 29 & 60 \\
\hline
\end{tabular}
\caption{Size and depth of certificates}
\label{tab3}
\end{table}

\begin{theorem}\label{th41}
A lower bound certificate $\sC$ exists for ones-ratio
$\alpha = \frac{\compnum}{\compden} \approx \comprat$.
\end{theorem}

\begin{proof}
The  certificate was found by computer search,
of size indicated in Table~\ref{tab3}. The computer
search required examining trees to depth $\compdep$, with
some paths having $\compones$ ones, although the path with
the largest ones ratio occured in a tree of depth $\comptd$,
having $\comptw$ ones. It is an interesting feature of this
certificate that the ``worst'' tree is not one of
maximal depth in the search.
\end{proof}

\noindent
{\em Proof of Theorem \ref{th11}.}
Since $a=1$ is in a cycle, we cannot apply Theorem \ref{th31}(i)
directly. Instead we
consider $n=41$ 
and note that $\rho (n) > \frac{55}{100}$
and $\gamma (n) > 20$.
Now we apply Theorem \ref{th31}(i) to $n=41$ using 
$\alpha =\frac{\compnum}{\compden}$ and the certificate of Theorem \ref{th41}.
Then the elements $n \to 41$ produced in Theorem \ref{th31}(i) have
$n \to 41 \to 1$ hence
$$\rho (n) \ge \min [\rho (n; 41), \rho (41)] \ge \frac{\compnum}{\compden},$$
as required. \hfill $\openbox$

The corresponding results for strong certificates are
given in Table~\ref{tab4}, giving the maximal value of $\alpha$
for which a strong certificate exists having trees with
all critical paths having at most $l$ ones. This table
is slightly less extensive because
larger searches were required. 
The breakpoint values of $\alpha$ are not the same as
for certificates. In many cases they do not
involve a larger depth for the largest tree, but sometimes
do.

\begin{table}
\renewcommand{\arraystretch}{0.625}
\begin{tabular}{|r|r|r|r|}
\hline
\rule[-0.07in]{0cm}{.23in}
ones-ratio $\alpha$ & size $|{\sC}|$ & max-weight $l$ & max-depth $k$ \\
\hline
$1/6   \approx 0.167$ &         6 &  1 &  6 \\
$1/4   =       0.250$ &        14 &  2 &  8 \\
$3/10  =       0.300$ &        26 &  3 & 10 \\
$4/13  \approx 0.308$ &        34 &  4 & 13 \\
$1/3   \approx 0.333$ &        36 &  5 & 15 \\
$6/17  \approx 0.353$ &        98 &  6 & 17 \\
$7/19  \approx 0.368$ &       204 &  7 & 19 \\
$8/21  \approx 0.381$ &       390 &  8 & 21 \\
$9/23  \approx 0.391$ &       848 &  9 & 23 \\
$2/5   =       0.400$ &       914 & 10 & 25 \\
$11/27 \approx 0.407$ &      2242 & 11 & 27 \\
$12/29 \approx 0.414$ &      4720 & 12 & 29 \\
$13/31 \approx 0.419$ &      8020 & 13 & 31 \\
$14/33 \approx 0.424$ &     16044 & 14 & 33 \\
$3/7   \approx 0.429$ &     16182 & 15 & 35 \\
$16/37 \approx 0.432$ &     34264 & 16 & 37 \\
$17/39 \approx 0.436$ &     64960 & 17 & 39 \\
$18/41 \approx 0.439$ &     91170 & 18 & 41 \\
$4/9   \approx 0.444$ &    158182 & 19 & 42 \\
$5/11  \approx 0.455$ &    838262 & 20 & 44 \\
$21/46 \approx 0.457$ &   1475962 & 21 & 46 \\
$11/24 \approx 0.458$ &   2374052 & 22 & 48 \\
$23/50 =       0.460$ &   4114846 & 23 & 50 \\
$6/13  \approx 0.462$ &   4114922 & 24 & 52 \\
$23/49 \approx 0.469$ &  25328092 & 25 & 53 \\
$26/55 \approx 0.473$ &  47636512 & 26 & 55 \\
$9/19  \approx 0.474$ &  47658126 & 27 & 57 \\
$28/59 \approx 0.475$ &  72183824 & 28 & 59 \\
\hline
\end{tabular}
\caption{Size and depth of strong certificates}
\label{tab4}
\end{table}

Let $\pi_{\alpha}(x)$ count the number 
of $n$ with 
$1 < n \le x $ which have finite total stopping time and
 ones-ratio $\rho(n) \ge \alpha.$

\begin{theorem}\label{th42}
A strong certificate $\sC$ exists for ones-ratio
$\alpha = \frac{\compsnum}{\compsden} \approx \compsrat$.
Furthermore there is a positive constant $c$ such that
$$\pi_{\alpha}(x) > c x^{1/\compsdep}.$$
\end{theorem}

\begin{proof}
The proof is similar to that for Theorem~\ref{th41}.
The bound for the exponent $x^{1/\compsdep}$ comes from (\ref{Eq311}).
\end{proof}

We remark on the size of the certificates versus
the computation needed to find them.
Although the certificates are large, they can be verified without 
computing the entire tree $\T_{[l]}^\ast (\bc )$;
one simply backtracks along the given test path(s).
However to find a minimal depth certificate, 
it appears necessary to calculate much
of the tree $\T_{[l]}^\ast ( \bc )$.
Thus a much larger computation seems needed to find 
the certificate then to check it afterwards.
Note that Lagarias and Weiss \cite[Theorem 3.1]{lw92}
showed that the expected number of leaves in a
``pruned'' $3x+1$ tree of depth $d$ is
$\left( \frac{4}{3} \right)^d$. \\

We now briefly describe the search method to find
the maximal value of $\alpha$ for a given level $l$,
the largest number of ones allowed in a critical path in
a certificate. We initialize the search with $\alpha$
taken to be the maximal value at level $l - 1$,
and with the given certificate.
Given a current test value of $\alpha$, 
a rational number, we
proceed with the algorithm above, looking for a certificate
at level $l$. If one exists, we determine the critical
path in the certificate having the largest ones-ratio,
which gave a value $\alpha'$, which is our current champion.
Then we take as a new test value $\alpha = \alpha' + \frac{1}{10000}$
and search it to level $l$ for a certificate.
If we get a new certificate by level $l$, we continue. Otherwise,
if some tree remains
unclosed at level $l$, we halt the computation, and
the current value of $\alpha'$ is maximal at level $l$.
This search procedure works because the maximal value $\alpha$
must be a rational number with denominator no larger than 
the maximal depth $k \le \frac {l}{\alpha}$, which in our
search is always less than $100$. All such values fall inside the
Farey sequence of order $100$, and the members of this Farey sequence
differ by more than $\frac{1}{10000}$.

The straightforward method of computation to determine
the critical path for a given $3x+1$ tree is to compute the
entire tree. If the maximal level $l$ to be searched
is known in advance (as it is in the incremental search algorithm above)
then we can get a speedup by pruning all branches of the tree at
the point where they cannot contribute any
critical path, i.e.  even if extended to the full allowed depth 
$k = \lfloor \frac{l}{\alpha} \rfloor $ 
with all edge labels $1$ they never have ones-ratio exceeding $\alpha.$

To conclude the computational results, 
we give statistics on the structure of the certificates found, as
described by the number of 
open vectors remaining at 
a given level $l' \le l$ of the certificate search; 
this value counts the number
of prefixes in the certificate having length at least $l'$.
We present  data on these values for different values of $\alpha$.
\begin{figure}[hbtp]
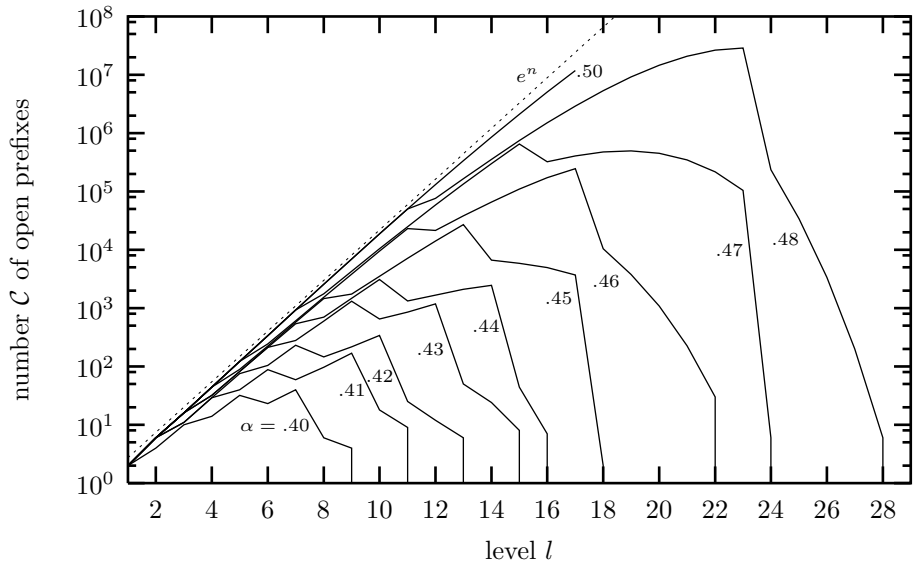
\centering
\input prefix_plot.tex
\caption{Number of open vectors in $\sC$ by depth, for various values
of $\alpha$}
\label{fig3}
\end{figure}
This data is plotted on a logarithmic scale in Figure~\ref{fig3}.
When extrapolated, this data makes some indication of what depth
of tree one is likely to have to search to get a certificate for
larger values of $\alpha$ than we obtained. The rate of increase
of $\alpha$ must slow down as $l$ increases, so we
view this extrapolation as indicating a {\em lower bound} on the 
depth to which one will have to search. We note, however, that the
data in Figure~\ref{tab3} exhibits some irregular
jumps  the value of $\alpha$, including a strange non-increase in
the value of $\alpha$ from level $l = 5$ to $l= 6$.

%
%

\section{Discussion}

(1) How large a lower bound can the method of this paper produce?

The stochastic models for $3x+1$ trees (branching random walk models)
studied in Lagarias and Weiss~\cite{lw92} predict that with
probability one a random tree grown to a large depth $k \to \infty$ 
will contain a path $\sP$ having $\gamma(\sP) > \gamma_{BP} - \epsilon,$
for any fixed $\epsilon > 0.$ However the method of proof
given above requires the lower bound $\gamma_0$ to be achieved for {\em all}
trees at some given depth~\footnote{
We actually require only that each such tree contain some subtree 
containing the root vertex, of possibly smaller depth,
having a path which achieves the given bound $\gamma \geq \gamma_0$}, and
not just for ``almost all'' trees of that depth. 
In terms of the stochastic models of \cite{lw92} and \cite{al95c}, 
we are sampling an
exponential number of trees (namely independently drawing
$3^k$ trees at depth $k$) out
of the double-exponential number of  possible trees in the 
branching random walk model at this depth,
and demanding that they all have
$\gamma(\sP) \geq \gamma_0$. In our deterministic situation, 
if this event ever occurs at some level $k$ , 
it yields a rigorous proof of the lower bound.
In the analogous stochastic model, we
should expect this cutoff event to occur {\em at some finite level $k$}
 with probability
one whenever the set of ``exceptional trees'' with 
$\gamma < \gamma_0$ has an exponentially small
probability, less than $3^{-k}$, for all sufficiently large $k$.
 In \cite{lw92} 
a large deviations analysis based on exactly this idea  produced the 
constant $\gamma_{BP}.$ We conclude that the stochastic
model supports the heuristic that this proof method should
be expected to work for any $\gamma_0 < \gamma_{BP}$, provided
that we can search trees to a sufficient depth, depending on
$\gamma_0.$  

To gain more confidence in this heuristic,
it would be interesting to make the stochastic model analysis of
\cite{lw92} more
precise, and to obtain a prediction
for each $\gamma_0$, of what depth $k$  one should expect to have to
go to before an independently
generated  sample of $3^k$ 
trees of that depth $k$ for the branching random walk,
 contains no maximal path with $\gamma(\sP) < \gamma_0$
with probability at least  $\frac {1}{2}.$

(2) Can a lower bound on the ones-ratio of $\frac{1}{2}$ be attained?

It remains surprising (to us) that a search to depth $\compdep$
of all trees was insufficient to produce a rigorous proof of 
the  bound 
$$\gamma_0 \geq \frac {2}{\log 4/3} \approx 6.95212, $$
a bound which is believed to apply to the vast majority of initial values $n$.
The statistical evidence provided by the cutoff values for increasing
values of the ones ratio from $0.41$ through $0.48$, when extrapolated,  
seem to indicate that a search to max-weight
$l$  of at least $38$, 
which corresponds to
trees of depth $76$, will be needed. 
We experimented on growing a few trees starting from the certificate
with ones-ratio $0.48$  to see how deep they had to become before
a path with ones ratio of $0.50$ or greater appeared; these seemed
to require depth around $68$. If this represents the ``knee'' of
the curve for $0.50$ in Figure~\ref{fig3}, then perhaps the depth necessary
for a certificate might be smaller, say $72$. In any case this 
would require a very large computation. 

With regard to extending the computations further, we note
that the certificate search algorithm of \S4 is well suited
to parallel and distributed computation, since each open
vector can be analyzed separately. Thus it should be possible
to search to considerably greater depth using a network
of machines. We leave a rigorous proof of the lower bound $\frac{1}{2}$
as a challenge to future researchers.

\ifx\undefined\bysame
\newcommand{\bysame}{\leavevmode\hbox to3em{\hrulefill}\,}
\fi

\end{document}